\documentclass[11pt]{article}

\usepackage[utf8]{inputenc}
\usepackage[T1]{fontenc}
\usepackage{mathpazo}
\usepackage{soul}
\RequirePackage[english]{babel}

\RequirePackage[a4paper,margin=2.5cm]{geometry}
\RequirePackage{parskip}

\newcommand{\affiliation}{\footnote}
\makeatletter
\def\@fnsymbol#1{\ensuremath{\ifcase#1\or *\or \dagger\or \ddagger\or \mathsection\or \|\or **\or \dagger\dagger \or \ddagger\ddagger \else\@ctrerr\fi}}
\makeatother

\usepackage{xcolor}
\definecolor{cblue}{RGB}{0,70,140}
\definecolor{cgreen}{RGB}{100,140,0}
\definecolor{cred}{RGB}{190,10,50}

\usepackage[shortlabels]{enumitem}
\setlist{itemsep=0ex,topsep=0ex,parsep=0.4ex}

\usepackage{amsmath,amsfonts,amssymb,amsthm,mathtools,bbm,centernot}

\usepackage[hyphens]{url}
\usepackage[hidelinks,colorlinks,pagebackref]{hyperref}
\hypersetup{citecolor=cgreen,linkcolor=cblue,urlcolor=cblue}
\usepackage[capitalise,nameinlink,noabbrev,compress]{cleveref}

\renewcommand*{\backref}[1]{}
\renewcommand*{\backrefalt}[4]{%
    \ifcase #1 Not cited.%
    \or $\uparrow$#2%
    \else $\uparrow$#2%
    \fi%
}

\let\oldbibliography\bibliography
\renewcommand{\bibliography}[1]{%
  {%
    \hypersetup{linkcolor=cred}%
    \bibliographystyle{bibstyle}%
    \oldbibliography{#1}%
  }%
}

\theoremstyle{plain}
\newtheorem{theorem}{Theorem}[section]
\newtheorem{lemma}[theorem]{Lemma}

\newtheorem{corollary}[theorem]{Corollary}

\newtheorem{remark}[theorem]{Remark}
\theoremstyle{definition}
\newtheorem{definition}[theorem]{Definition}
\newtheorem*{keylem}{Key Lemma}
\newtheorem*{genkeylem}{Generalized Key Lemma}
\makeatletter
\renewenvironment{proof}[1][\proofname]
{\par\pushQED{\qed}
    \normalfont\topsep6\p@\@plus6\p@\relax\trivlist
    \item[\hskip\labelsep\bfseries#1\@addpunct{.}]
    \ignorespaces}
{\popQED\endtrivlist\@endpefalse}
\makeatother

\let\emptyset\varnothing

\newcommand{\arb}{\operatorname{arb}}
\newcommand{\mad}{\operatorname{mad}}
\newcommand{\cpack}{\chi^{\star}_{c}}
\newcommand{\fpack}{\chi^{\bullet}_{c}}

\newcommand{\cstar}{\chi_c^{\star}}
\newcommand{\chip}{\chi_c^\star}
\newcommand{\chif}{\chi_c^\bullet}

\title{List packing for graphs with small arboricity}

\author{
S. Cambie\thanks{Department of Computer Science, KU Leuven Campus Kulak-Kortrijk, 8500 Kortrijk, Belgium. Supported by FWO grants 1225224N and 1215727N. Email: \protect\href{mailto:stijn.cambie@hotmail.com}{\textsf{stijn.cambie@hotmail.com}}.}
\and
W. Cames van Batenburg\affiliation{D\'epartement d'Informatique, Universit\'e libre de Bruxelles, Belgium (\textsf{\href{mailto:w.p.s.camesvanbatenburg@gmail.com}{w.p.s.camesvanbatenburg@gmail.com}}). Supported by the Belgian National Fund for Scientific Research (FNRS).}
\and
X. Zhu\affiliation{School of Mathematical Sciences, Zhejiang Normal University, China (\textsf{\href{mailto:xdzhu@zjnu.edu.cn}{xdzhu@zjnu.edu.cn}}). Supported by Natural Science Foundation of China grant:    NSFC 12371359.}
}
\date{\today}

\begin{document}

\maketitle

\begin{abstract}
We prove that every graph $G$ whose edge set is the union of $r$ forests has list packing number at most $2^r$.
This means that for all $k\ge 2^r$ and every list-assignment $L$ with lists of size  $k$, the graph $G$ admits $k$ disjoint proper $L$-colourings. 
While for 
large $r$
this 
bound is far from optimal, for $r\in\{1,2\}$ it is sharp.
As a direct corollary of the case $r=2$, we obtain that every triangle-free planar graph has list packing number at most $4$, thus improving upon a result of
Dvo\v{r}\'ak, Masa{\v{r}}{\'i}k, Mus{\'i}lek and Pangr{\'a}c
and providing the optimal answer to a question of Dvo\v{r}\'ak, Postle and Norin on (weighted) $\epsilon$-flexibility of such graphs.
As a corollary of the case $r=3$, we obtain a short alternative proof that every $K_5$-minor-free (and hence every planar) graph has list-packing number at most $8$, matching the state of the art. All results above also hold for the stronger notion of correspondence packing.
\end{abstract}

\section{Introduction}

We present a follow-up 
of an earlier paper~\cite{CambieCamesZhu23} by the same authors on disjoint list colourings in planar graphs. All graphs in this manuscript are simple and finite. A \emph{$k$-list-assignment} $L$ assigns to every vertex $v$ a
private set $L(v)$ of $k$ allowed colours, and an \emph{$L$-colouring} is a proper vertex-colouring $c:V(G) \rightarrow \bigcup_{v} L(v)$ such that $c(v)\in L(v)$ for every vertex $v$. The smallest $k$ such that every $k$-list-assignment of $G$ admits an $L$-colouring is the \emph{list chromatic number} $\chi_{\ell}(G)$. A particularly useful application of this widely studied parameter has been the short inductive proof~\cite{Thomassen94} that planar graphs have list chromatic number (and hence chromatic number) at most $5$. Likewise, this note will find its applications in planar graphs.

In the study of graph colourings, more often than not, one is already satisfied if \emph{one}  colouring with the desired properties is proved to exist. We are instead interested in finding multiple colourings that are moreover pairwise disjoint. In this vein, \emph{list packings} were introduced in~\cite{CambieCamesDaviesKang24}, as an interpolation between list colourings, independent transversals and strong colourings.
More precisely, an \emph{$L$-packing} of $G$ is a
collection of $k$ pairwise disjoint $L$-colourings $c_1,c_2,\ldots, c_k$, meaning that $c_i(v)\neq c_j(v)$ for every distinct $i,j \in [k]$ and $v \in V(G)$.
As a consequence, these $k$ colourings partition every list $L(v)$. The \emph{list packing number} $\chi_{\ell}^{\star}(G)$ of $G$ is the smallest integer $k$ such that every $k$-list-assignment $L$ of $G$ admits an $L$-packing.
 The most basic but still completely open problem in this area is whether $\chi_{\ell}^{\star}(G) \le C \cdot \chi_{\ell}(G)$ for some absolute constant $C>0$.

Beyond the notion of list colourings, we will work in the more general correspondence setting of Dvo\v{r}\'ak
and Postle \cite{DvorakPostle18}.  A \emph{correspondence-cover} of
$G$ is a pair $\mathcal H=(L,H)$, where $H$ is a graph, the sets
$\{L(v):v\in V(G)\}$ partition $V(H)$, every $H[L(v)]$ is a
clique, and, for each edge $uv\in E(G)$, the edges of $H$ between
$L(u)$ and $L(v)$ form a matching $M_{uv}$. 
If $uv\notin E(G)$, then there is no edge between $L(u)$ and $L(v)$.  
The cover is \emph{$k$-fold} if
$\lvert L(v)\rvert=k$ for every $v$.  An
$\mathcal H$-colouring is an independent transversal of $\mathcal{H}$, meaning an independent set of $H$ that intersects every list $L(v)$. A $\mathcal H$-packing is a collection of $k$ pairwise disjoint
independent transversals.  Thus, an $\mathcal H$-packing partitions
$V(H)$ into $k$ independent transversals.  The
\emph{correspondence packing number} $\cpack(G)$ is the least
integer $k$ such that every $k$-fold correspondence-cover of
$G$ has a packing.  Every list-assignment has an associated
correspondence-cover, so 
$\chi_{\ell}(G)\le \chi_{\ell}^{\star}(G) \le \chi_c^{\star}(G)$ for every graph $G$.


Since their introduction, the packing numbers have been studied in terms of various graph invariants, including maximum degree, degeneracy, fractional chromatic number, treewidth and pathwidth. For more background and related papers we refer to~\cite{CambieCamesDaviesKang24,CCvBDK2,CambieCames26,CambieCamesZhu23,CamesCranstonKardos26,KMZ26,KMMP22,KM25,CranstonSmithRoberge25,CT25}. Here we study yet another parameter, which turns out to be particularly useful for applications to planar graphs: the \emph{arboricity} $\arb(G)$ of a graph $G$ is the least number of forests whose edge
sets partition $E(G)$.  Our main result gives a direct connection
between this classical density parameter and correspondence packing.

\begin{theorem}\label{thm:arboricity}
Every graph $G$ satisfies
$$
  \cpack(G)\leq 2^{\arb(G)}.
$$
\end{theorem}

As we will see later, Theorem~\ref{thm:arboricity} is sharp for $\arb(G) \in \{1,2\}$. For $\arb(G)=3$, the optimal bound could be $6,7$, or $8$. For $\arb(G)\ge 4$, a simple degeneracy argument already supplies something better. Indeed,
the degeneracy of $G$ is at most $2\arb(G)-1$ and as proved in~\cite{CambieCamesDaviesKang24}, twice degeneracy is a sharp upper bound on $\chi_c^{\star}(G)$.
Despite the bound of Theorem~\ref{thm:arboricity} being very bad for large arboricity, we have kept the statement in this general form because it reflects the proof method, which could be interpreted as a careful product colouring over the forests. Moreover, this form might inspire the use of other decomposition theorems beyond forests. For the nontrivial applications of Theorem~\ref{thm:arboricity}, we are primarily interested in graphs with arboricity $2$ or $3$.\\

Recall that the \emph{maximum average degree}
$\mad(G)
  =
  \max_{\varnothing\ne J\subseteq G}
  \frac{2\lvert E(J)\rvert}{\lvert V(J)\rvert}
$
is defined as the average degree maximized over all subgraphs of $G$. The Nash--Williams arboricity theorem~\cite{NashWilliams64}  relates arboricity with maximum average degree, turning Theorem~\ref{thm:arboricity} into a density bound

\begin{corollary}\label{cor:madsummary}
$$ \chi_c^{\star}(G) \leq \begin{cases}
 8 & \text{ if } \mad(G)<5+\frac{3}{7}, \\
 4 & \text{ if } \mad(G)<3+\frac{3}{5},\\
 2 & \text{ if } \mad(G)<2.
\end{cases}$$
\end{corollary}

The bound $8$ was already known under the weaker condition $\mad(G)<6$, but with a long proof that required a computer-assisted preprocessing step to count certain tuples of small permutations~\cite{CambieCamesZhu23}. 
Our condition $\mad(G)<3+\frac{3}{5}$ is sharp due to the complete graph on $5$ vertices with one edge removed,
and improves the previous best-known~\cite{CambieCamesZhu23, CranstonSmithRoberge25} condition $\mad(G)<10/3$. The result under condition $\mad(G)<2$ is trivial, but it allows us to emphasize a curious gap~\cite{CH23}: there exists no graph $G$ with $\chi_c^{\star}(G)=3$. For completeness, we also mention the result from~\cite{CambieCamesZhu23, CranstonSmithRoberge25} that $\chi_c^{\star}(G)\le 5$ whenever $\mad(G)<4$, which is sharp because $\chi_c^{\star}(K_5)=6$.
Our second corollary recovers results from~\cite{CambieCamesZhu23} for graphs that exclude a small clique minor. 

\begin{corollary}\label{thr:minorfreesummary}
$$ \chi_c^{\star}(G) \leq \begin{cases}
 8 & \text{ if } G \text{ is } K_5\text{-minor-free,}\\
 4 &\text{ if } G \text{ is } K_4\text{-minor-free,}\\
 2 & \text{ if } G \text{ is } K_3\text{-minor-free}.
\end{cases}$$
\end{corollary}

The bounds for $K_4$-minor-free and $K_3$-minor-free graphs are optimal. Recent work~\cite{CamesCranstonKardos26} shows that
there exist  $K_5$-minor-free graphs with correspondence packing number $6$, and every $6$-fold correspondence-cover of every $K_5$-minor-free graph $G$ admits three
pairwise disjoint colourings, but $\chi_c^{\star}(G)\le 6$
remains conjectural. Our third corollary is perhaps the most appealing.
\begin{corollary}\label{cor:planar-minor}
 $$ \chi_c^{\star}(G) \leq \begin{cases}
  8 & \text{ if } G \text{ is planar,}\\
 4 & \text{ if } G \text{ is triangle-free planar,}\\
2 & \text{ if } G \text{ is a forest}.
\end{cases}
$$
\end{corollary}

For triangle-free planar graphs, Corollary~\ref{cor:planar-minor} improves the previous upper bound \cite{CambieCamesZhu23,CranstonSmithRoberge25} from $5$ to $4$. This was open even for bipartite planar graphs and it is sharp in a strong sense. For instance, every cycle has correspondence packing number $4$, and due to a construction in~\cite{CambieCamesZhu23}, there exist planar graphs of arbitrarily large girth with list packing number $4$. \\
The result for general planar graphs recovers the
known bound $8$ for planar graphs
\cite{CambieCamesZhu23,CranstonSmithRoberge25}, but via a shorter proof that avoids extensive case analysis of matchings in auxiliary bipartite graphs.

There is also an optimal consequence for \emph{flexibility}, a notion introduced by Dvo\v{r}\'ak, Norin and Postle~ \cite{DvorakNorinPostle19} which measures to what extent one can force an $L$-colouring to agree with predetermined preferred colours.  Let $L$ be a
list-assignment of $G$.  A \emph{request} is a map $r$ with domain
$D\subseteq V(G)$ such that $r(v)\in L(v)$ for every $v\in D$.
The graph $G$ is \emph{$\varepsilon$-flexible} with respect to $L$ if
some $L$-colouring satisfies at least $\varepsilon\lvert D\rvert$
requests, for every $r$. 
Moreover, a graph $G$ is \emph{weighted $\epsilon$-flexible} with respect to $L$ if there exists some probability distribution on the $L$-colourings of $G$ such that for every vertex $v$, and every colour $x$ in its list $L(v)$, the vertex $v$ is assigned colour $x$ with probability at least $\epsilon$. Clearly $\epsilon$ cannot be larger than 
$\frac{1}{|L(v)|}$,
since otherwise the expected value of the number of colours assigned to $v$ would be larger than $1$.
Moreover, weighted $\varepsilon$-flexibility implies
$\varepsilon$-flexibility~\cite{DvorakNorinPostle19}, as can be seen by computing the expected value of the number of satisfied requests.\\

Given a list packing for a $k$-list-assignment $L$ of a graph $G$, the uniform probability distribution on the $k$ disjoint $L$-colourings shows that $G$ is weighted $\frac{1}{k}$-flexible with respect to $L$. Therefore Corollary~\ref{cor:planar-minor} implies the next result. 

\begin{corollary}\label{cor:flexibility}
Every triangle-free planar graph is weighted $1/4$-flexible, and
hence $1/4$-flexible, with respect to every $4$-list-assignment.
The constant $1/4$ is best possible in the weighted statement.
\end{corollary}

In their foundational paper, Dvo\v{r}\'ak, Norin and Postle asked~ \cite[Problem~5]{DvorakNorinPostle19}
whether there exists an absolute constant $\varepsilon>0$ such that
every triangle-free planar graph is weighted $\epsilon$-flexible for some $\epsilon>0$ with respect to every $4$-list-assignment. Dvo\v{r}\'ak,
Masa\v{r}\'{\i}k, Mus\'{\i}lek and Pangr\'ac answered this
affirmatively \cite{DvorakMasarikMusilekPangrac21} via an involved discharging argument and reducible configurations. The parameters in
their proof give $\varepsilon=2^{-186}$.  Our result strengthens
their theorem to the largest possible value $\varepsilon=1/4$.
More substantially, the four colourings may be chosen to form a
packing, and the underlying packing theorem extends from list colouring to arbitrary
correspondence colouring.
For more results on (weighted) flexibility that were obtained independent of the packing perspective, particularly for planar graphs, we refer to~\cite{BCK25,DMMP20, CCFHMM22, BMS22, BB23,DvorakMasarikMusilekPangrac21}.

In Section~\ref{sec:main-proof}
we prove Theorem~\ref{thm:arboricity} via a Key Lemma.  The corollaries are proved in Section~\ref{sec:corollaries}. In Section~\ref{sec:generalizedkeylemma}, we present a Generalized Key Lemma that goes beyond forest decompositions. We expect it to find future applications to more involved graph classes. We finish
with a few open problems in Section~\ref{sec:concluding}.

\section{Proof of the main theorem}\label{sec:main-proof}
The proof of the main theorem relies on the Key Lemma below. 
For a spanning subgraph $J\subseteq G$, write $G-E(J)$ for the spanning graph with edge set $E(G)\setminus E(J)$.

\begin{keylem}\label{lem:forest-extension}
Let $J$ be a spanning subgraph of $G$ such that $G-E(J)$ is a forest.  Then
$$\cstar(G)\le 2\cstar(J)$$
\end{keylem}

\begin{proof}
We first argue that the lists can be split in some suitable way.  Fix an integer $k\ge1$ and a $2k$-fold correspondence cover $\mathcal H=(L,H)$ of $G$.  By exploiting that $F:=G-E(J)$ is a $1$-degenerate graph, we shall partition every list $L(v)$ as the disjoint union of two lists $L_0(v)$ and  $L_1(v)$ of equal size $|L_0(v)|=|L_1(v)|=k$,
in such a way that for every edge $uv$ of $F$, every edge of $H$ between $L(u)$ and $L(v)$ is either between $L_0(u)$ and $L_1(v)$ or between $L_1(u)$ and $L_0(v)$. In other words, for every $uv\in E(F)$, we want that no edge of the matching between $L(u)$ and $L(v)$ has both endpoints in $L_0(u) \cup L_0(v)$ or both endpoints in $L_1(u)\cup L_1(v)$.


Root every component of $F$.  At a root $v$, choose an arbitrary equipartition of $L(v)$.  Suppose $u$ is the parent of $v$ and the partition of $L(u)$ is already fixed.  If $x\in L(v)$ is matched, along the edge $uv$, to $y\in L_i(u)$, then force $x$ into $L_{1-i}(v)$.  Since the edges between $L(u)$ and $L(v)$ form a matching, at most $k$ colours of $L(v)$ are forced into either part.  Assign the remaining colours (if there are any) so that both parts have size $k$.  Proceeding away from the roots gives the desired split of the lists.

Let $\mathcal H_i=\mathcal H[\bigcup_v L_i(v)]$ be the cover restricted to part $i$. Since no edge of $\mathcal{H}$ arising from $F$ has both ends in part $i$, it follows that $\mathcal H_i$ is a $k$-fold cover of $J$.  Moreover, every independent transversal of $\mathcal H_i$ is automatically an independent transversal of $\mathcal{H}$.

Now set $k=\cstar(J)$.  Each $\mathcal H_i$ has a packing into $k$ independent transversals, and the packings of $\mathcal{H}_1$ and $\mathcal{H}_2$ together form a packing of the original $2k$-fold cover of $G$.  Hence $\cstar(G)\le2k$.

\end{proof}

\begin{corollary}\label{cor:iterate}
If $J\subseteq G$ is spanning and $G-E(J)$ is the union of $r$ forests, then $\cstar(G)\le 2^r \cdot \cstar(J)$.
\end{corollary}

\begin{proof}
Add the forests one at a time and iteratively apply the Key Lemma.
\end{proof}

\begin{proof}[Proof of Theorem~\ref{thm:arboricity}]
Take $J$ edge-less in Corollary~\ref{cor:iterate}.
\end{proof}

\section{Proofs of the corollaries}\label{sec:corollaries}

\begin{proof}[Proof of Corollary~\ref{cor:madsummary}]
By the Nash--Williams arboricity theorem~\cite{NashWilliams64}, $\arb(G)\leq t$ if and
only if every subgraph $J$ with at least two vertices satisfies
$
  |E(J)|\leq t\cdot (|V(J)|-1).
$
If $\arb(G)\geq4$, then some subgraph $J$ with
$n=|V(J)|$ and $m=|E(J)|$ has $m > 3(n-1)$.  This is only possible if $n\geq7$, so
$
  \mad(G)\geq\frac{2m}{n}\geq6-\frac4n\geq\frac{38}{7}.
$
Thus $\mad(G)<38/7$ implies $\arb(G)\leq3$.\\
Similarly, if
$\arb(G)\geq3$, then some subgraph has $m\geq2n-1$. Then $n\geq5$, so
$
  \mad(G)\geq\frac{2m}{n}\geq4-\frac2n\geq\frac{18}{5}.
$
Thus $\mad(G)<18/5$ implies $\arb(G)\leq2$.  Finally, $\mad(G)<2$ implies $G$ is a forest and thus $\arb(G)\leq1$. The three bounds
now follow from Theorem~\ref{thm:arboricity}.
\end{proof}

\begin{proof}[Proof of Corollary~\ref{thr:minorfreesummary}]
 The extremal bounds of Mader~\cite{Mader68} give
$|E(J)|\leq3|V(J)|-6$ for every $K_5$-minor-free subgraph $J$
with at least three vertices, and $|E(J)|\leq2|V(J)|-3$ for every
$K_4$-minor-free subgraph $J$ with at least two vertices. 
Nash--Williams therefore gives
arboricity at most $3$ and $2$,
respectively.  A $K_3$-minor-free graph is a forest, and hence has
arboricity at most $1$.  Apply Theorem~\ref{thm:arboricity}.
\end{proof}

\begin{proof}[Proof of Corollary~\ref{cor:planar-minor}]
Every planar graph is $K_5$-minor-free, so the first assertion follows
from Corollary~\ref{thr:minorfreesummary}.  If $G$ is triangle-free
and planar, then every subgraph $J$ with $n\geq3$ vertices has
$|E(J)|\leq2n-4\leq2(n-1)$ by Euler's formula. The same bound holds for $n=2$. Hence $\arb(G)\leq2$ by Nash--Williams. 
\end{proof}

\begin{proof}[Proof of Corollary~\ref{cor:flexibility}]
Let $L$ be a $4$-list-assignment of a triangle-free planar graph
$G$.  By Corollary~\ref{cor:planar-minor}, the associated
correspondence-cover has a packing of four $L$-colourings.  At every
vertex these colourings use each colour in its list exactly once, so
the uniform distribution on them assigns every colour probability
$1/4$.  Thus $G$ is weighted $1/4$-flexible, and hence
$1/4$-flexible.  Conversely, for any probability distribution on
$L$-colourings and any vertex $v$, the probabilities of the four
colours in $L(v)$ sum to $1$. Therefore no larger uniform lower
bound is possible.
\end{proof}

\section{Generalized key lemma}\label{sec:generalizedkeylemma}
The results in this section are not necessary for the main results of the current paper, but we present them with a view to potential future applications to more involved graph classes, via partitions into graphs other than forests. Since forests are $1$-degenerate graphs, our Key Lemma is a special case of Corollary~\ref{cor:degenerate} below. 
First we need some definitions.


\begin{definition}
The \emph{fractional list packing number} $\chi_{\ell}^{\bullet}(G)$ of a graph $G$ is the smallest integer $k$ such that $G$ is weighted $\frac{1}{k}$-flexible with respect to every $k$-fold list-assignment.
\end{definition}
Similar to how the notion of weighted $\frac{1}{k}$-flexibility can be viewed as a relaxation of a list packing, one can relax a correspondence packing to a fractional correspondence packing.
\begin{definition}
The \emph{fractional correspondence packing number} $\fpack(G)$ of a graph $G$ is the smallest integer $k$ such that every $k$-fold cover $(L,H)$ of $G$ admits a probability distribution on the independent transversals such that for every vertex $v$, every colour in $L(v)$
appears with probability $1/k$.
\end{definition}

For background and results on the fractional packing numbers, we refer to~\cite{CCvBDK2,CambieCamesZhu23, CambieCames26}. The next definition is new, as far as we are aware. It facilitates a streamlined statement of the Generalized Key Lemma.
\begin{definition}
We define the \emph{correspondence splitting number} $\sigma_c(F)$ of a graph $F$ as the least positive
integer $s$ such that the following holds for every integer
$k\ge 1$. For every $sk$-fold correspondence cover $\mathcal H=(L,H)$ of $F$, there are partitions
$L(v)=L_1(v)\cup \cdots \cup L_s(v)$ of the lists, such that for every vertex $v$ and every $j\in[s]$ we have  $|L_j(v)|=k$, and for every edge $uv\in E(F)$ and every $j\in[s]$, there is
no edge in $H$ between $L_j(u)$ and $L_j(v)$.\\
The \emph{list splitting number} $\sigma_{\ell}(F)$ of a graph $F$ is defined analogously; then the conditions only are required to hold for the correspondence covers that arise from an $sk$-list assignment.
\end{definition}

It is easy to see that $\sigma_{\ell}(G)\le \sigma_{c}(G)$ and  $\chi_{\ell}^{\bullet}(G)\le \chi_c^{\bullet}(G)\leq \chi_c^{\star}(G)$ and $\chi_{\ell}(G)\le \chi_{\ell}^{\bullet}(G)\leq \chi_{\ell}^{\star}(G) \le \chi_{c}^{\star}(G) \le \sigma_c(G)$ for every graph $G$. 

\begin{genkeylem}\label{lem:genkeylem}
Let $J$ be a spanning subgraph of $G$. Write $F:=G-E(J)$. Then
$$\chip(G)\le \sigma_c(F)\cdot\chip(J),$$
$$\chif(G)\le \sigma_c(F)\cdot \chif(J),$$
$$\chi_{\ell}^{\star}(G)\le \sigma_{\ell}(F)\cdot\chi_{\ell}^{\star}(J),$$
$$
\chi_{\ell}^{\bullet}(G)\le \sigma_{\ell}(F)\cdot \chi_{\ell}^{\bullet}(J).
$$
\end{genkeylem}

\begin{proof}
We will focus on the first inequality, as the other three are proved similarly. Let $s:=\sigma_c(F)$ and $k:=\chip(J).$
Take an arbitrary $sk$-fold correspondence cover $(L,H)$ of $G$. By
definition of $\sigma_c(F)$, there exists a partition of $V(H)$ into parts $P_1,P_2,\ldots, P_s$, such that for every $j\in [s]$ the following holds: 
for every vertex $v$ the truncated list $L_j(v):= L(v)\cap P_j$ has size $k$, and for every edge $uv \in E(F)$ there is no edge between $L_j(u)$ and $L_j(v)$ in $H$.

For each $j\in [s]$, the induced subgraph $H_j:=H[P_j]$ equipped with lists $L_j:=(L_j(v))_{v\in V(G)}$ forms a $k$-fold correspondence
cover $(L_j,H_j)$ of $J$ and thus admits a packing, a partition into $k$ independent
transversals. As no edge of $F$ corresponds to any edge in $H_j$,
these are also independent transversals of $(L,H)$. Taking all $s$ packings together yields
an $sk$-packing of $(L,H)$. Hence $\chip(G)\le s \cdot k.$

The fractional statement in terms of $\chi_c^{\bullet}$ follows identically, by first choosing $j\in [s]$ uniformly at random and then sampling a fractional packing from the cover $(L_j,H_j)$ of $J$. In the second stage (conditional on the value of $j$) every vertex of $H_j$ is selected with probability $1/k$ because $\chi_c^{\bullet}(J)=k$. Every vertex of $H$ belongs to
exactly one $H_j$, and hence altogether is selected with probability
$\frac1{s}\cdot\frac1k=\frac1{sk}$.

The last two inequalities in terms of $\chi_{\ell}^{\star}$ and $\chi_{\ell}^{\bullet}$ also follow via the same proof, by restricting to the correspondence covers that arise from list-assignments.
\end{proof}

The next lemma allows to convert the Generalized Key Lemma into a concrete special case, allowing decompositions into graphs with bounded degeneracy. By definition, a graph is \emph{$d$-degenerate} if there exists an ordering of the vertices such that every vertex has at most $d$ neighbours earlier in the ordering.

\begin{lemma}\label{lem:degeneracy_splitting}
For every $d$-degenerate graph $F$, we have $\sigma_c(F)\le 2d$. 
\end{lemma}
\begin{proof}

Fix an integer $k\ge 1$, and let $\mathcal H=(L,H)$ be an arbitrary
$2dk$-fold correspondence cover of $F$. We need to show that for every
$v\in V(F)$ one can partition the list
$L(v)$ into $2d$ disjoint sublists $L_1(v),L_2(v),\ldots,L_{2d}(v)$, each of size $k$,
in such a way that for every $uv\in E(F)$ and every $j\in[2d]$,
there is no edge between $L_j(u)$ and $L_j(v)$ in $H$.

Choose a degeneracy ordering $v_1,\ldots,v_n$ of $F$, so that every $v_s$ has at most $d$ neighbours among $v_1,\ldots,v_{s-1}$. We construct the partitions in this order.
Suppose that $v=v_s$ is being processed, and let $u_1,\ldots,u_r$
be its earlier neighbours in $F$. Their partitions have already
been fixed. We now need to construct $L_1(v),\ldots, L_{2d}(v)$.

For a colour $x\in L(v)$ and an index $j\in[2d]$, say that $j$ is
\emph{forbidden for $x$} 
if $x$ has a neighbour in $\bigcup_{i\in [r]} L_j(u_i)$, in the cover graph $H$. Otherwise $j$ is \emph{allowed} for $x$.

Since the edges between $L(u_i)$ and $L(v)$ form
a matching, each earlier neighbour $u_i$ forbids at most one index
$j$ for a fixed colour $x$. Hence every $x\in L(v)$ has at most
$d$ forbidden indices, and consequently at least $d$ allowed
indices.
Conversely, fix $j\in[2d]$. For each earlier neighbour $u_i$,
the sublist $L_j(u_i)$ contains $k$ colours, so at most $k$ colours
of $L(v)$ are matched to colours in $L_j(u_i)$. Thus at most
$ rk\le dk $
colours of $L(v)$ are forbidden to appear in $L_j(v)$.

We now use Hall's theorem. Construct an auxiliary bipartite graph $B$ whose
left part is the list of colours $L(v)$, and whose right part consists of $k$ copies
$j^1,\ldots,j^k$
of every index $j\in[2d]$. Join $x\in L(v)$ to every copy of $j$
precisely when $j$ is allowed for $x$. Note that both sides of $B$ have size $2dk$.
We verify Hall's condition.
Let $S\subseteq L(v)$. If
$|S|\le dk$,
then every $x\in S$ is adjacent to all $k$ copies of at least
$d$ indices. Hence
$|N_B(S)|\ge dk\ge |S|.$
Now suppose instead that $|S|>dk$.
Fix $j\in[2d]$. At most $dk$ colours of $L(v)$ are forbidden
to appear in $L_j$, so some colour of $S$ is allowed in $L_j$. Hence all $k$ copies $j^1,\ldots,j^k$
are neighbours of $S$. Since this holds for every $j\in[2d]$, it follows that $N_B(S)$
is the entire right side of $B$, so $|N_B(S)|=2dk\ge |S|.$
Thus we may apply Hall's matching theorem: $B$ has a perfect matching.

We convert the perfect matching to a choice of sublists of $L(v)$. Assign colour $x\in L(v)$ to sublist $L_j(v)$ if and only if $x$ is matched in $B$ to one of the $k$ copies of $j$. Since $j$ is allowed for $x$, this does not create any undesirable edge between $L_j(u)$ and $L_j(v)$ for any edge $uv\in E(F)$.
Moreover, since each of the $k$ copies of $j$ is matched, this gives
$|L_j(v)|=k$
for every $j$.
\end{proof}

\begin{corollary}
\label{cor:degenerate}
Let $J$ be a spanning subgraph of a graph $G$ such that $F:=G-E(J)$ is $d$-degenerate for some $d\ge 1$. Then $\chip(G)\le 2d\cdot\chip(J)$ and $\chif(G)\le 2d\cdot \chif(J)$.
\end{corollary}
\begin{proof}
Apply the Generalized Key Lemma to the bound in Lemma~\ref{lem:degeneracy_splitting}.  
\end{proof}

\begin{remark}\label{rem:wouldbeniceproduct}\qquad
\begin{itemize}
    \item  Taking $J$ the spanning edgeless graph in Corollary~\ref{cor:degenerate}, we  recover the result from~\cite{CambieCamesDaviesKang24} that $\chi_c^{\star}(G)$ is at most twice the degeneracy of $G$. This is sharp due to large complete bipartite graphs.

\item Taking $J$ a spanning forest (which is $1$-degenerate) in Corollary~\ref{cor:degenerate}, we  recover the Key Lemma.

\item Taking $k=1$ in the definition of $\sigma_c(F)$ yields
$\chip(F)\le \sigma_c(F)$.
It would be interesting to understand under which conditions there is equality, $\chip(F)=\sigma_c(F)$.
Let $\mathcal{F}$ denote the class of such graphs. Then for every spanning subgraph $J$ of $G$ such that $G-E(J)\in \mathcal{F}$, the following attractive product inequality holds: 
$$\chip(G)\le \chip(G-E(J)) \cdot \chip(J).$$
\end{itemize}
\end{remark}

We end this section with a refined variant of the Generalized Key Lemma. Given a list assignment $L$ and positive integer $k$, a \emph{$k$-fold $L$-colouring} is a mapping $\phi$ such that $\phi(v)$ is a size $k$-subset of $L(v)$ for each vertex $v$, and $\phi(u) \cap \phi(v) = \emptyset$ for each edge $uv$. For a positive integer $s$, we say a graph G is \emph{list $(sk,k)$-packable} 
if for any $sk$-fold assignment $L$, there are  $k$-fold $L$-colourings $\phi_1, \ldots, \phi_s$   such that $\{\phi_1(v),\ldots, \phi_s(v)\}$ is a partition of $L(v)$ for each vertex $v$.  The concept of correspondence-$(sk,k)$-packable is defined similarly. 
By the proof of the Generalized Key Lemma, it follows that if a graph $G$ decomposes into $G_1 \cup G_2$, and $G_1$ is $(sk,k)$-packable and $G_2$ is $k$-packable, then $G$ is $sk$-packable.

One might even hope that $G_1$ $s$-packable and $G_2$ $k$-packable implies that $G$ is $sk$-packable. To achieve that, a natural conjecture would be that  if $G$ is $s$-packable, then $G$ is $(sk,k)$-packable for any positive integer $k$. However, that is probably too optimistic because a closely related statement for list colourings is false: Erd\H{o}s, Rubin and Taylor conjectured~\cite{ERT80} that if $G$ is $(a,b)$-choosable, then $G$ is $(ak,bk)$-choosable for any positive integer $k$, and this has been disproved~\cite{DHS19}.

\section{Concluding remarks}\label{sec:concluding}

We conclude with open problems. Some are achievable with assistance from current AI  models. However, the challenge is to obtain an elegant humanly readable proof.
  
\begin{enumerate}[label=\normalfont(\arabic*)]
\item What is the optimal upper bound on $\chi_c^*(G)$ for graphs $G$ with arboricity $3$? Is it $6,7$ or $8$?
  \item
   From~\cite{CambieCamesZhu23}, it is known that $\fpack(G) \le 3$
   for planar graphs of girth $\ge 6$, while there are triangle-free planar graphs with $\fpack(G)= 4$.  Do all planar graphs of girth $5$ satisfy
  $\fpack(G)\leq3$?
  For the integral packing numbers this is false: there exist planar graphs $G$ of arbitrarily large girth with  $\chi_{\ell}^{\star}(G) = \chi_{c}^{\star}(G) = 4$.

  \item In line with Remark~\ref{rem:wouldbeniceproduct}, what other graph decompositions lead to good upper bounds on the various packing numbers? 
  Is there any graph $F$ for which $\sigma_c(F)\neq \chi_c^{\star}(F)$, or $\sigma_{\ell}(F)\neq \chi_{\ell}^{\star}(F)$? 

  \item Regarding planar graphs, the main remaining open problem is to determine the optimal upper bound on $\chi_{\ell}^{\star}(G)$ or $\chi_{c}^{\star}(G)$ without any girth constraints. 
  In upcoming work~\cite{C26},  a planar graph with correspondence packing number $6$ is presented. There are now three independent proofs that $\chi_{c}^{\star}(G)\le 8$ for all planar graphs, the other two proofs appearing in~\cite{CambieCamesZhu23,CranstonSmithRoberge25}. 
  An obstruction for the decompositions presented in the current paper is that they involve (iterated applications of) bounds of the form $ \sigma \cdot  \chi_c^{\star}(J)$ for some integer $\sigma \ge 2$ and spanning subgraph $J$. Since the only graphs with correspondence packing number $<4$ are forests, this cannot achieve anything better than $2\cdot 4=8$ for planar graphs. 
  \item Graphs with maximum average degree $<6$ have~\cite{CambieCamesZhu23} correspondence packing number $\le 8$. Not every such graph can be edge-partitioned into at most three forests, so the current approach fails to recover that result in an easy way. It is however interesting to note that every such graph can be partitioned into at most one forest and two pseudo-forests. Here, a pseudo-forest is a graph such that every connected component has at most one cycle. Can this lead to an \emph{elegant} alternative proof?
\item For $k\in \{2,4,5,8\}$, we now have a reasonably good picture of the optimal conditions on maximum average degree that ensure $\chi_c^{\star}(G)\le k$. For $k\in\{2,4,5\}$ the bounds are sharp.
What about $k\in \{6,7\}$? 
We reiterate a question from~\cite{CambieCamesZhu23}: does $\mad(G)<k-1$ imply $\chi_c^{\star}(G)\le k$? 
We conjecture that this is true for $k=6.$
\end{enumerate}

\paragraph{Disclosure and background}
In conversation with a closed-weight AI model about proving $\chi_c^{\star}(G)\le 4$ for every triangle-free planar graph $G$,
 one of the suggestions of the model was to use that they admit a decomposition into two forests.
This eventually led to a direct but obscure proof by induction for triangle-free planar graphs, based on binary sequences associated with the vertices. The authors isolated the key step, formulated the Key Lemma, and from there obtained the other corollaries and the Generalized Key Lemma.


\bibliography{bib}
\end{document}